\documentclass[12pt,twoside]{article}

\usepackage[frame,cmtip,arrow,matrix,line,graph,curve]{xy}
\usepackage{amssymb}

\date{}

\begin{document}

\centerline{\bf }

\centerline{\bf }

\centerline{}

\centerline{}

\centerline{\Large{\bf Anti-invariant Riemannian submersions }}

\centerline{}

\centerline{\Large{\bf  from locally conformal Kaehler manifolds }}

\centerline{}

\centerline{\bf {Majid Ali Choudhary}}

\centerline{}

\centerline{Department of Mathematics}

\centerline{Zakir Husain Delhi College (E), Delhi University, New Delhi -110002(India)   }

\centerline{E-mail : majid\_alichoudhary@yahoo.co.in}

\centerline{}

\newtheorem{Theorem}{\quad Theorem}[section]

\newtheorem{Definition}[Theorem]{\quad Definition}

\newtheorem{Corollary}[Theorem]{\quad Corollary}

\newtheorem{Lemma}[Theorem]{\quad Lemma}

\newtheorem{Example}[Theorem]{\quad Example}

\centerline{}

\begin{abstract}  B. Sahin \cite{11} introduced the notion of anti-invariant Riemannian submersions from almost Hermitian manifolds onto Riemannian manifolds. In the present paper we extend the notion of anti-invariant and Lagrangian Riemannian submersions (a special anti-invariant Riemannian submersion) to the case of locally conformal Kaehler manifolds. We discuss the geometry of foliation and obtain some decomposition theorems for the total manifold of such submersions.
\end{abstract}

{\bf MSC}  53C15, 53B20, 53C43. \\

{\bf Keywords:} Anti-invariant Riemannian submersions, Lagrangian Riemannian submersions, locally conformal Kaehler manifold.
\section{Introduction}
\setcounter{equation}{0}
\renewcommand{\theequation}{1.\arabic{equation}}
\hspace {1cm}
    Locally conformal Kaehler manifolds (Shortly, l.c.K. manifolds) have been rich source of research for many geometers. A number of papers have appeared on these manifolds and their submanifolds (See: \cite{13}, \cite{6} for details). On the other hand, the study of the Riemannian submersions $\pi:M\rightarrow B$ of a Riemannian manifold $M$ onto a Riemannian manifold $B$ was initiated by B. O'Neil \cite{4}. Gray \cite{101}, Ianus \cite{102}, Park (\cite{103}, \cite{104}), Sahin (\cite{105}, \cite{106}) etc. have also investigated submersions in different settings. In 2010, B. Sahin \cite{11} introduced the notion of anti-invariant Riemannian submersions from almost Hermitian manifolds onto Riemannian manifolds. Here, in the present paper we extend the notion of anti-invariant and Lagrangian Riemannian submersions (a special anti-invariant Riemannian submersion) to the case when the total manifold is locally conformal Kaehler manifold. We obtain some integrability conditions for the horizontal distribution while it is noted that vertical distribution is always integrable. We also discuss necessary and sufficient condition for a Lagrangian Riemannian submersions to be totally geodesic and obtain decomposition theorems for the total manifold of such submersions.

\section{Preliminaries}
\setcounter{equation}{0}
\renewcommand{\theequation}{2.\arabic{equation}}

In this section, first we define l.c.K. manifold, recall the notion of Riemannian submersion between Riemannian manifolds. Then, we give brief review of basic facts of basic vector fields.
\begin{Definition} \cite{6} Let $(\widetilde{M},J,g)$  be a Hermitian manifold of dimension $2m$. Let $\Omega$ be the Kaehler $2$-form associated with $J$ and $g$ i.e. $\Omega(X,Y)=g(X,JY)$  for all $X,Y \in \chi(\widetilde{M})$. The manifold $\widetilde{M}$  is called locally conformal Kaehler (l.c.K.) manifold if there is a closed $1$-form $\omega$ defined globally on $\widetilde{M}$  such that $d\Omega=\omega\wedge\Omega.$ \end{Definition}
The closed $1$-form $\omega$  is called the Lee form of the l.c.K. manifold $\widetilde{M}$. Also $(\widetilde{M},J,g)$ is globally conformal Kahler (g.c.K.) (respectively Kahler) if the Lee form $\omega$ is exact (respectively $\omega = 0$). Note that any simply connected l.c.K is g.c.K. For a l.c.K. manifold $(\widetilde{M},J,g)$ we define the Lee vector field $B_1 =\omega^\sharp$ where $\sharp$ means the rising of the indices with respect to $g$, namely $g(X,B_1)=\omega(X)$; for all $X \in \chi(\widetilde{M})$. If $\widetilde{\nabla}$ denotes the Levi Civita connection of $(\widetilde{M},J,g)$ then we have
\begin{eqnarray}(\widetilde{\nabla}_X J)Y=\frac{1}{2} \{\theta(Y)X-\omega(Y)JX-g(X,Y)A-\Omega(X,Y)B_1\}\end{eqnarray}
for any  $X,Y \in \chi(\widetilde{M})$. Here, $\theta=\omega oJ$  and  $A=-JB_1$ are the anti-Lee form and the anti-Lee vector field, respectively \cite{6}.

\centerline{}

Let $(M^{m},g)$ and $(B^b,g_B)$ be Riemannian manifolds, where $dim (M)=m$, $dim (B)=b$ and $m>b$. A Riemannian submersion $\pi: M\rightarrow B$ is a map of $M$ onto $B$ satisfying the following axioms:\begin{description}
                                           \item[(1)] $\pi$ has maximal rank.
                                           \item[(2)] The differential $\pi_*$ preserves the lengths of horizontal vectors.
                                         \end{description}
For each $q \in B, \pi^{-1}(q)$ is an $(m-b)$ dimensional submanifold of $M$. The submanifold $\pi^{-1}(q), q \in B$ is called fiber. A vector field on $M$ is called vertical if it is always tangent to fibers. A vector field on $M$ is called horizontal if it is always orthogonal to fibers. In order to relate the geometry of $M$ and $B$ we have the concept of basic vector fields defined as follows:
\begin{Definition} A vector field $X$ on $M$ is said to be basic provided $X$ is horizontal and is $\pi$-related to a vector field $X_*$ on $B$, that is, $\pi_{*}oX=X_{*}o\pi$.  \end{Definition} It is to be noted that we denote the projection morphisms on the distributions $ker\pi_*$ and $(ker\pi_*)^\bot$ by $\mathcal{V}$ and $\mathcal{H}$, respectively. We have the following lemma for basic vector fields \cite{4}
\begin{Lemma} Let $\pi: M\rightarrow B$ be a Riemannian submersion between Riemannian manifolds $M$ and $B$ and $X,Y$ be basic vector fields on $M$, then
\begin{description}
        \item[(a)] $g(X,Y)=g_B(X_*,Y_*)o\pi$.
        \item[(b)] $\mathcal{H}[X,Y]$ of $[X,Y]$ is a basic vector field and corresponds to $[X_*,Y_*]$, that is, $([X,Y]^{\mathcal{H}})=(X_*,Y_*).$
        \item[(c)] $[V,X]$ is vertical for any vertical vector $V$.
        \item[(d)] $\mathcal{H}(\nabla_X Y)$ is the basic vector field corresponding to $\nabla^{*}_{X_{*}} Y_{*}$, where $\nabla^{*}$ is a Levi-Civita connection on $B$.
      \end{description}
 \end{Lemma}
The geometry of Riemannian submersions is characterized by O'Neills tensors $\mathcal{T}$ and $\mathcal{A}$ defined for vector fields $E,F$ on $M$ by \cite{4}

\begin{eqnarray} \mathcal{A}_E F=\mathcal{H}\nabla_{\mathcal{H}E}\mathcal{V}F+\mathcal{V}\nabla_{\mathcal{H}E}\mathcal{H}F \end{eqnarray}
\begin{eqnarray} \mathcal{T}_E F=\mathcal{H}\nabla_{\mathcal{V}E}\mathcal{V}F+\mathcal{V}\nabla_{\mathcal{V}E}\mathcal{H}F \end{eqnarray}
where $\nabla$ is the Levi-Civita connection of $g$. It is easy to see that a Riemannian submersion $\pi:M\rightarrow B$ has totally geodesic fibres if and only if $\mathcal{T}$ vanishes identically. For any $E \in \Gamma(TM), \mathcal{T}_E$ and $\mathcal{A}_E$ are skew-symmetric operators on $(\Gamma(TM),g)$ reversing the horizontal and vertical distributions. It is also easy to see that $\mathcal{T}$ is vertical, $\mathcal{T}_E =\mathcal{T}_{\mathcal{V}E}$ and $\mathcal{A}$ is horizontal, $\mathcal{A}=\mathcal{A}_{\mathcal{H}E}$. We note that tensor fields $\mathcal{T}$ and $\mathcal{A}$ satisfy (\cite{4}, \cite{11})

\begin{eqnarray} \mathcal{T}_U W=\mathcal{T}_W U, \forall U,W \in \Gamma(ker\pi_*) \end{eqnarray}
\begin{eqnarray} \mathcal{A}_X Y=-\mathcal{A}_Y X=\frac{1}{2}\mathcal{V}[X,Y], \forall X,Y \in (\Gamma(ker\pi_*)^\bot). \end{eqnarray}

Now we state the following lemma \cite{11}
\begin{Lemma} For $X, Y \in (\Gamma(ker\pi_*)^\bot)$ and $V,W \in \Gamma(ker\pi_*)$, we have the following relations:
\begin{description}
        \item[(a)] $\nabla_V W=\mathcal{T}_V W +\hat{\nabla}_V W$
        \item[(b)] $\nabla_V X=\mathcal{H}\nabla_V X+\mathcal{T}_V X$
        \item[(c)] $\nabla_X V=\mathcal{A}_X V+\mathcal{V}\nabla_X V$
        \item[(d)] $\nabla_X Y=\mathcal{H}\nabla_X Y+\mathcal{A}_X Y$      \end{description}
where $\hat{\nabla}_V W=\mathcal{V}\nabla_V W$. If $X$ is basic, then $\mathcal{H}\nabla_V X=\mathcal{A}_X V$.
 \end{Lemma}

\section{Anti-invariant Riemannian submersions}
\setcounter{equation}{0}
\renewcommand{\theequation}{3.\arabic{equation}}

In this section, we recall the definitions of anti-invariant Riemannian submersion and Lagrangian Riemannian submersion, investigate the integrability of distributions and obtain a necessary and sufficient condition for such submersion to be totally geodesic map. We also investigate the harmonicity of a Lagrangian Riemannian submersion.

B. Sahin \cite{11} defined anti-invariant Riemannian submersion by the following way
\begin{Definition} Let $M$ be a complex $m$-dimensional almost Hermitian manifold with Hermitian metric $g$ and almost complex structure $J$ and $B$ be a Riemannian manifold with Riemannian metric $g_B$. Suppose that there exists a Riemannian submersion $\pi:M\rightarrow B$ such that $ker\pi_*$ is anti-invariant with respect to $J$, that is $J(ker\pi_*)\subseteq(ker\pi_*)^\bot$. Then we say that $\pi$ is an anti-invariant Riemannian submersion. \end{Definition}
Let $\pi:(M,g,J)\rightarrow (B,g_B)$ be an anti-invariant Riemannian submersion from an almost Hermitian manifold $M$ onto a Riemannian manifold $B$. From the above definition we observe that $J(ker\pi_*)^\bot\cap(ker\pi_*)\neq{0}$ and hence we have
\begin{eqnarray}(ker\pi_*)^\bot=J(ker\pi_*)\oplus\mu \end{eqnarray}
where $\mu$ denotes the orhogonal complementary distribution to $J(ker\pi_*)$ in $(ker\pi_*)^\bot$ and it is invariant under $J$. Thus for any $X \in \Gamma((ker\pi_*)^\bot)$ we have
\begin{eqnarray}JX=BX+CX\end{eqnarray}
where $BX \in \Gamma(ker\pi_*)$ and $CX \in \Gamma(\mu)$. On the other hand, since we have $\pi_*((ker\pi_*)^\bot)=TB$ and $\pi$ is a Riemannian submersion using (3.2) it can be shown that $g_B(\pi_* JV,\pi_* CX)=0$ for any $X \in \Gamma((ker\pi_*)^\bot)$ and $V \in \Gamma(ker\pi_*)$, which implies that
\begin{eqnarray}TB=\pi_*(J(ker\pi_*))\oplus\pi_*(\mu).\end{eqnarray}

Finally, we recall the notion of harmonic maps between Riemannian manifolds \cite{3}. Let us suppose that $(M,g)$ and $(B,g_B)$ be Riemannian manifolds and  $\phi:M\rightarrow B$ be a smooth map. Then the differential $\phi_*$ of $\phi$ can be viewed as a section of the bundle Hom$(TM,\phi^{-1}(TB))\rightarrow M$, where $\phi^{-1}(TB)$ is the pullback bundle which has fibres $(\phi^{-1}(TB))_p=T_{\phi(p)} B, p \in M$. $Hom(TM,\phi^{-1}(TB))$ has a connection $\nabla$ induced from Levi-Civita connection $\nabla^M$ and the pullback connection. The second fundamental form of $\phi$ is given by
\begin{eqnarray}(\nabla\phi_*)(X,Y)=\nabla^\phi_X \phi_*(Y)-\phi_*(\nabla^M_X Y) \end{eqnarray} for
$X, Y \in \Gamma(TM)$, where $\nabla^\phi$ is the pullback connection.

We give the following lemma to be used later in order to prove the theorems.

\begin{Lemma}
Let $\pi$ be an anti-invariant Riemannian submersion from a locally conformal Kaehler manifold $(M,g,J)$ to a Riemannian manifold $(B,g_B)$. Then we have
\begin{description}
  \item[(i)] $g(CY,JV)=0$
  \item[(ii)] $g(\nabla_X CY,JV)=-g(CY,J\mathcal{A}_X V)+\frac{1}{2}\omega(V)g(CY,CX)$
  \item[(iii)] $g(\nabla_V BY,CX)=g(CX,\mathcal{T}_V BY)=-g(BY,\mathcal{T}_V CX)$
\end{description}
where $\omega$ is closed $1$-form defined globally on $M$ and for $X,Y \in \Gamma((ker\pi_*)^\bot), V \in \Gamma(ker\pi_*).$
\end{Lemma}

\noindent

{\bf Proof} (i) For any $Y \in \Gamma((ker\pi_*)^\bot)$ and $V \in \Gamma(ker\pi_*)$ from (3.2) we have
\begin{eqnarray*}
g(CY,JV)&=&g(JY-BY,JV)\\
&=&g(JY,JV)
\end{eqnarray*}
as $BY \in \Gamma(ker\pi_*)$ and $JV \in \Gamma((ker\pi_*)^\bot)$. But, $ g(JY,JV)=g(Y,V)=0$ and hence the result.

(ii) For any $X, Y \in \Gamma((ker\pi_*)^\bot)$ and $V \in \Gamma(ker\pi_*)$ using part (i) and equation (2.1) we have
\begin{eqnarray*}
g(\nabla_X CY,JV)&=&-g(CY,\nabla_X JV)\\
&=&-g(CY,J\nabla_X V)+\frac{1}{2}\omega(V)g(CY,JX)
\end{eqnarray*}
where we have assumed that $B_1 \in \Gamma(ker\pi_*)$. Now using (3.2) we have
\begin{eqnarray*}
g(\nabla_X CY,JV)&=&-g(CY,J\nabla_X V)+\frac{1}{2}\omega(V)g(CY,BX+CX)\\
&=&-g(CY,J\nabla_X V)+\frac{1}{2}\omega(V)g(CY,CX)
\end{eqnarray*}
as $CY \in \Gamma(\mu)$ and $BX \in \Gamma(ker\pi_*)$. Applying Lemma 2.4 we have
\begin{eqnarray*}
g(\nabla_X CY,JV)=-g(CY,J\mathcal{A}_X V)-(CY,J\mathcal{V}\nabla_X V)+\frac{1}{2}\omega(V)g(CY,CX)
\end{eqnarray*}
since $J\mathcal{V}\nabla_X V \in \Gamma(Jker\pi_*)$, we obtain
\begin{eqnarray*}
g(\nabla_X CY,JV)=-g(CY,J\mathcal{A}_X V)+\frac{1}{2}\omega(V)g(CY,CX)
\end{eqnarray*}
which proves the result.
\hfill$\Box$

\centerline{}
It is to be noted that from now we will assume that $B_1 \in (ker\pi_*)$ through out the paper.
During the proofs whenever it seems necessary, we have supposed horizontal vector fields to be basic. Now, we discuss the integrability of the horizontal distribution $(ker\pi_*)^\bot$. One can notice that distribution $ker\pi_*$ is integrable.

\begin{Theorem}
Let $(M,g,J)$ be a locally conformal Kaehler manifold and $\pi$ be an anti-invariant Riemannian submersion from $M$ onto a Riemannian manifold $(B,g_B)$. Then the following assertions are equivalent to each other:
\begin{description}
  \item[(a)] $(ker\pi_*)^\bot$ is integrable
  \item[(b)] $g_B((\nabla\pi_*)(Y,BX),\pi_* JV)=g_B((\nabla\pi_*)(X,BY),\pi_* JV)+g(CY,J\mathcal{A}_X V) \\
~~~~~~~~~~~~~~~~~~~~~~~~~~~~~~~~-g(CX,J\mathcal{A}_Y V)-\frac{1}{2}g(BY,B_1)g(X,JV)\\
~~~~~~~~~~~~~~~~~~~~~~~~~~~~~~~~+\frac{1}{2}g(BX,B_1)g(Y,JV)$
  \item[(c)] $g(\mathcal{A}_Y BX-\mathcal{A}_X BY,JV)=-g(CY,J\mathcal{A}_X V)+g(CX,J\mathcal{A}_Y V)\\
~~~~~~~~~~~~~~~~~~~~~~~~~~~~~~~~+\frac{1}{2}g(BY,B_1)g(X,JV)-\frac{1}{2}g(BX,B_1)g(Y,JV)$
\end{description}
for $X,Y \in \Gamma((ker\pi_*)^\bot)$ and $V \in \Gamma(ker\pi_*).$
\end{Theorem}
{\bf Proof.} If $Y \in \Gamma((ker\pi_*)^\bot)$ and $V \in \Gamma(ker\pi_*)$, then from definition 3.1 it follows that $JY \in \Gamma(ker\pi_*\oplus\mu)$ and $JV \in \Gamma((ker\pi_*)^\bot)$. Thus using (2.1) for $X \in \Gamma((ker\pi_*)^\bot)$ we have
\begin{eqnarray*}
g([X,Y],V)&=&g(J[X,Y],JV)\\
&=&g(J\nabla_X Y,JV)-g(J\nabla_Y X,JV)\\
&=&g(\nabla_X JY,JV)-\frac{1}{2}\theta(Y)g(X,JV)\\
&&-g(\nabla_Y JX,JV)+\frac{1}{2}\theta(X)g(Y,JV)
\end{eqnarray*}
here $\theta=\omega oJ$,  $\Omega(X,Y)=g(X,JY)$ and $g(X,B_1)=\omega(X)$, then the assumption that $B_1 \in \Gamma(ker\pi_*)$ gives
\begin{eqnarray*}
g([X,Y],V)&=&g(\nabla_X JY,JV)-g(\nabla_Y JX,JV)-\frac{1}{2}g(BY,B_1)g(X,JV)\\
&&+\frac{1}{2}g(BX,B_1)g(Y,JV)\\
&=&g(\nabla_X BY,JV)+g(\nabla_X CY,JV)-g(\nabla_Y BX,JV)\\
&&-g(\nabla_Y CX,JV)-\frac{1}{2}g(BY,B_1)g(X,JV)+\frac{1}{2}g(BX,B_1)g(Y,JV)\\
\end{eqnarray*}
where we have used (3.2). As $\pi$ is a Riemannian submersion, we obtain
\begin{eqnarray*}
g([X,Y],V)&=&g(\pi_*\nabla_X BY,\pi_*JV)+g(\nabla_X CY,JV)
-g_B(\pi_*\nabla_Y BX,\pi_*JV)\\
&&-g(\nabla_Y CX,JV)-\frac{1}{2}g(BY,B_1)g(X,JV)+\frac{1}{2}g(BX,B_1)g(Y,JV)
\end{eqnarray*}
Taking into account the definition of second fundamental form and using lemma 3.2, we have
\begin{eqnarray*}
g([X,Y],V)&=&g_B(-(\nabla\pi_*)(X,BY)+(\nabla\pi_*)(Y,BX),\pi_*JV)\\
&&-g(CY,J\mathcal{A}_X V)+g(CX,J\mathcal{A}_Y V)\\
&&-\frac{1}{2}g(BY,B_1)g(X,JV)+\frac{1}{2}g(BX,B_1)g(Y,JV)
\end{eqnarray*}
Thus, we have shown that (a)$\Leftrightarrow$(b).

Next, for $X,Y \in \Gamma((ker\pi_*)^\bot)$ and $V \in \Gamma(ker\pi_*)$ taking into consideration Lemma 2.4 , we have
\begin{eqnarray*}
(\nabla\pi_*)(X,BY)&-&(\nabla\pi_*)(Y,BX)\\
&=&-\pi_*(\nabla_X BY)+\pi_*(\nabla_Y BX)\\
&=&-\pi_*(\nabla_X BY-\nabla_Y BX)\\
&=&\pi_*(\mathcal{A}_Y BX-\mathcal{A}_X BY)
\end{eqnarray*}
Hence, we get
\begin{eqnarray*}
g_B((\nabla\pi_*)(X,BY)&-&(\nabla\pi_*)(Y,BX),\pi_* JV)\\
&=&g_B(\pi_*(\mathcal{A}_Y BX-\mathcal{A}_X BY),\pi_* JV)\\
&=&g(\mathcal{A}_Y BX-\mathcal{A}_X BY,JV)
\end{eqnarray*}
as $\mathcal{A}_X BY-\mathcal{A}_Y BX \in \Gamma((ker\pi_*)^\bot)$, we conclude that (b)$\Leftrightarrow$(c).
\hfill$\Box$

\begin{Definition}\cite {11} An anti-invariant Riemannian submersion $\pi$ is said to be Lagrangian Riemannian submersion if $J(ker\pi_*)=(ker\pi_*)^\bot$. If $\mu\neq \{0\}$, then $\pi$ is said to be proper anti-invariant Riemannian submersion.
\end{Definition}

From Theorem 3.3, we have the following result for Lagrangian Riemannian submersion

\begin{Corollary}
Let $(M,g,J)$ be a locally conformal Kaehler manifold and $\pi:M\rightarrow B$ be Lagrangian Riemannian submersion of $M$ onto a Riemannian manifold $(B,g_B)$. Then the following are equivalent:
\begin{description}
  \item[(a)] $(ker\pi_*)^\bot$ is integrable
  \item[(b)] $(\nabla\pi_*)(X,JY)=(\nabla\pi_*)(Y,JX)-\frac{1}{2}g(BY,B_1)X+\frac{1}{2}g(BX,B_1)Y$
  \item[(c)] $\pi_*(\mathcal{A}_X JY-\mathcal{A}_Y JX)=\frac{1}{2}g(BY,B_1)X-\frac{1}{2}g(BX,B_1)Y$
\end{description}
for $X,Y \in \Gamma((ker\pi_*)^\bot)$.
\end{Corollary}

{\bf Proof.}
               For any $X,Y \in \Gamma((ker\pi_*)^\bot)$ and $V \in \Gamma(ker\pi_*)$ we observe that $JX \in \Gamma(ker\pi_*)$ and $JV \in \Gamma((ker\pi_*)^\bot)$. Then in the light of (2.1) we see
\begin{eqnarray*}
g([X,Y],V)&=&g(J[X,Y],JV)\\
&=&g(J\nabla_X Y,JV)-g(J\nabla_Y X,JV)\\
&=&g(\nabla_X JY,JV)-g(\nabla_Y JX,JV)\\
&&-\frac{1}{2}g(BY,B_1)g(X,JV)+\frac{1}{2}g(BX,B_1)g(Y,JV)
\end{eqnarray*}
Taking into account that $\pi$ is a Riemannian submersion and using (3.4), we arrive at
\begin{eqnarray*}
g([X,Y],V)&=&g_B(\pi_* \nabla_X JY,\pi_* JV)-g_B(\pi_* \nabla_Y JX,\pi_* JV)\\
&&-\frac{1}{2}g(BY,B_1)g(X,JV)+\frac{1}{2}g(BX,B_1)g(Y,JV)\\
&=&-g_B((\nabla \pi_*)(X, JY),\pi_* JV)+g_B((\nabla \pi_*)(Y, JX),\pi_* JV)\\
&&-\frac{1}{2}g(BY,B_1)g(X,JV)+\frac{1}{2}g(BX,B_1)g(Y,JV)
\end{eqnarray*}
Hence, $(ker\pi_*)^\bot$ is integrable if and only if
\begin{eqnarray*}
g_B((\nabla \pi_*)(X, JY),\pi_* JV)&=&g_B((\nabla \pi_*)(Y, JX),\pi_* JV)-\frac{1}{2}g(BY,B_1)g(X,JV)\\
&&+\frac{1}{2}g(BX,B_1)g(Y,JV)
\end{eqnarray*}
This shows that(a)$\Leftrightarrow$(b).

Now, in view of (3.4) for $X,Y \in \Gamma((ker\pi_*)^\bot)$ we have
\begin{eqnarray*}
(\nabla\pi_*)(Y,JX)&-&(\nabla\pi_*)(X,JY)\\
&=&-\pi_*(\nabla_Y JX)+\pi_*(\nabla_X JY)\\
&=&\pi_*(\mathcal{H}(\nabla_X JY)-\mathcal{H}(\nabla_Y JX))\\
&=&\pi_*(\mathcal{A}_X JY-\mathcal{A}_Y JX)
\end{eqnarray*}
                Hence, we conclude that (b)$\Leftrightarrow$(c).
\hfill$\Box$

\section{Totally geodesicness on $M$}
\setcounter{equation}{0}
\renewcommand{\theequation}{4.\arabic{equation}}

In this section we study geometry of leaves of $(ker\pi_*)$ and $(ker\pi_*)^\bot$ of anti-invariant Riemannian submersion and Lagrangian Riemannian submersion. First, we have the following result for the totally geodesicness of $(ker\pi_*)^\bot$.

\begin{Theorem}
Let $(M,g,J)$ be a locally conformal Kaehler manifold and $\pi$ be an anti-invariant Riemannian submersion from $M$ onto a Riemannian manifold $(B,g_B)$. Then the following assertions are equivalent to each other:
\begin{description}
  \item[(a)] $(ker\pi_*)^\bot$ defines a totally geodesic foliation on $M$
  \item[(b)] $g(\mathcal{A}_X BY,JV)=g(CY,J\mathcal{A}_X V)-\frac{1}{2}\omega(V)g(CY,CX)\\
~~~~~~~~~~~~~~~~~~~~ +\frac{1}{2}g(BY,B_1)g(X,JV)+\frac{1}{2}g(X,Y)g(B_1,V)$
  \item[(c)] $g_B((\nabla\pi_*)(X,JY),\pi_* JV)=-g(CY,J\mathcal{A}_X V)+\frac{1}{2}\omega(V)g(CY,CX)\\
~~~~~~~~~~~~~~~~~~~~~~~~~~~~~~~~ -\frac{1}{2}g(BY,B_1)g(X,JV)-\frac{1}{2}g(X,Y)g(B_1,V)$
\end{description}
for $X,Y \in \Gamma((ker\pi_*)^\bot)$ and $V \in \Gamma(ker\pi_*).$
\end{Theorem}

{\bf Proof.}
                   In view of (2.1) and (3.2) for $X,Y \in \Gamma((ker\pi_*)^\bot)$ and $V \in \Gamma(ker\pi_*)$, we have
\begin{eqnarray*}
g(\nabla_X Y,V)&=&g(J\nabla_X Y,JV)\\
&=&g(\nabla_X JY,JV)-\frac{1}{2}g(BY,B_1)g(X,JV)-\frac{1}{2}g(X,Y)g(B_1,V)\\
&=&g(\nabla_X BY,JV)+g(\nabla_X CY,JV)-\frac{1}{2}g(BY,B_1)g(X,JV)\\
&&-\frac{1}{2}g(X,Y)g(B_1,V)\\
&=&g(\mathcal{A}_X BY,JV)-g(CY,J\mathcal{A}_X V)+\frac{1}{2}\omega(V)g(CY,CX)\\
&&-\frac{1}{2}g(BY,B_1)g(X,JV)-\frac{1}{2}g(X,Y)g(B_1,V)
\end{eqnarray*}
where we have used Lemma 2.4 and Lemma 3.2. Therefore, $(ker\pi_*)^\bot$ defines a totally geodesic foliation on $M$ if and only if
\begin{eqnarray*}
g(\mathcal{A}_X BY,JV)&=&g(CY,J\mathcal{A}_X V)-\frac{1}{2}\omega(V)g(CY,CX)\\
&&+\frac{1}{2}g(BY,B_1)g(X,JV)+\frac{1}{2}g(X,Y)g(B_1,V)
\end{eqnarray*}
Thus, we conclude that (a)$\Leftrightarrow$(b).
Next, in the light of equation (3.4), we have
\begin{eqnarray*}
g(\mathcal{A}_X BY,JV)&=&g(\nabla_X BY,JV)\\
&=&g(\nabla_X JY,JV)-g(\nabla_X CY,JV)\\
&=&g_B(\pi_*\nabla_X JY,\pi_*JV)-g(\nabla_X CY,JV)\\
&=&-g_B((\nabla\pi_*)(X,JY),\pi_*JV)+g_B(\nabla^{\pi}_X \pi_*(JY),\pi_* JV)\\
&&-g(\nabla_X CY,JV)\\
&=&-g_B((\nabla\pi_*)(X,JY),\pi_*JV)+g(\nabla_X CY,JV)-g(\nabla_X CY,JV)\\
&=&-g_B((\nabla\pi_*)(X,JY),\pi_*JV)
\end{eqnarray*}
hence, we conclude that (b)$\Leftrightarrow$(c).
\hfill$\Box$

If $\pi$ is a Lagrangian Riemannian submersion, then we have the following corollary giving characterization for totally geodesicness of horizontal distribution of Lagrangian Riemannian submersion.

\begin{Corollary}
Let $\pi$ be a Lagrangian Riemannian submersion from a locally conformal Kaehler manifold $(M,g,J)$ onto a Riemannian manifold $(B,g_B)$. Then the following assertions are equivalent to each other:
\begin{description}
  \item[(a)] $(ker\pi_*)^\bot$ defines a totally geodesic foliation on $M$
  \item[(b)] $g_B(\mathcal{A}_X JY,JV)=\frac{1}{2}g(JY,B_1)g(X,JV)+\frac{1}{2}g(X,Y)g(B_1,V)$
  \item[(c)] $g_B((\nabla\pi_*)(X,JY),\pi_* JV)=-\frac{1}{2}g(JY,B_1)g(X,JV)-\frac{1}{2}g(X,Y)g(B_1,V)$
\end{description}
for $X,Y \in \Gamma((ker\pi_*)^\bot)$ and $V \in \Gamma(ker\pi_*).$
\end{Corollary}

{\bf Proof.}
                In view of (2.1) for $X,Y \in \Gamma((ker\pi_*)^\bot)$ and $V \in \Gamma(ker\pi_*)$, we have
\begin{eqnarray*}
g(\nabla_X Y,V)&=&g(J\nabla_X Y,JV)\\
&=&g(\nabla_X JY,JV)-\frac{1}{2}\theta(Y)g(X,JV)-\frac{1}{2}g(X,Y)g(B_1,V)\\
&=&g_B(\pi_* \nabla_X JY,\pi_* JV)-\frac{1}{2}\theta(Y)g(X,JV)-\frac{1}{2}g(X,Y)g(B_1,V)\\
&=&g_B(\pi_* (\mathcal{A}_X JY),\pi_* JV)-\frac{1}{2}\theta(Y)g(X,JV)-\frac{1}{2}g(X,Y)g(B_1,V)
\end{eqnarray*}
 Therefore, $(ker\pi_*)^\bot$ defines a totally geodesic foliation on $M$ if and only if
\begin{eqnarray*}g_B(\mathcal{A}_X JY,JV)=\frac{1}{2}\theta(Y)g(X,JV)+\frac{1}{2}g(X,Y)g(B_1,V)\end{eqnarray*}
Thus, (a)$\Leftrightarrow$(b).
Next, taking into account (3.4) we have
\begin{eqnarray*}g_B(\mathcal{A}_X JY,JV)&=&g_B(\nabla_X JY,JV)\\
&=&g_B(\pi_* \nabla_X JY,\pi_* JV)\\
&=&-g_B((\nabla\pi_*)(X,JY),\pi_* JV)
\end{eqnarray*}
which shows that
\begin{eqnarray*}g_B((\nabla\pi_*)(X,JY),\pi_* JV)=-\frac{1}{2}\theta(Y)g(X,JV)-\frac{1}{2}g(X,Y)g(B_1,V)
\end{eqnarray*}
proving that (b)$\Leftrightarrow$(c).
 \hfill$\Box$

Now, let us suppose that $X \in \Gamma(\mu)$ and $V \in \Gamma(ker\pi_*)$. Then, by the definition of second fundamental form (3.4), we have
\begin{eqnarray*}
(\nabla\pi_*)(V,X)=\nabla^\pi_V \pi_* X-\pi_* \nabla_V X.
\end{eqnarray*}
Further , we have
\begin{eqnarray*}
(\nabla\pi_*)(X,V)=\nabla^\pi_X \pi_* V-\pi_* \nabla_X V.
\end{eqnarray*}
Since, second fundamental form is symmetric, from above two equations we conclude that
\begin{eqnarray}
\nabla^\pi_V \pi_* X=0.
\end{eqnarray}
For the totally geodesicness of the leaves of $(ker\pi_*)$, we have

\begin{Theorem}
Let $(M,g,J)$ be a locally conformal Kaehler manifold and $\pi$ be an anti-invariant Riemannian submersion from $M$ onto a Riemannian manifold $(B,g_B)$. Then the following assertions are equivalent to each other:
\begin{description}
  \item[(a)] $(ker\pi_*)$ defines a totally geodesic foliation on $M$
  \item[(b)] $\mathcal{T}_V BX+\mathcal{A}_CX V=0$ or $\mathcal{T}_V BX+\mathcal{A}_CX V \in \Gamma(\mu)$
  \item[(c)] $g_B((\nabla\pi_*)(V,JX),\pi_* JW)=0$
\end{description}
for $X \in \Gamma((ker\pi_*)^\bot)$ and $V,W \in \Gamma(ker\pi_*).$
\end{Theorem}

{\bf Proof.} For any $X \in \Gamma((ker\pi_*)^\bot)$ and $V,W \in \Gamma(ker\pi_*)$, we have
\begin{eqnarray*}
g(\nabla_V W,X)&=&g(J\nabla_V W,JX)\\
&=&g(\nabla_V JW,JX)\\
&=&-g(JW,\nabla_V JX)
\end{eqnarray*}
where we have used (2.1) and the fact that $(ker\pi_*)$ and $(ker\pi_*)^\bot$ are orthogonal. Now, taking into consideration equation (3.2) and applying Lemma 2.4, we observe that
\begin{eqnarray*}
g(\nabla_V W,X)&=&-g(JW,\nabla_V BX)-g(JW,\nabla_V CX)\\
&=&-g(JW,\mathcal{T}_V BX)-g(JW,\mathcal{A}_{CX} V)\\
&=&-g(JW,\mathcal{T}_V BX+\mathcal{A}_{CX} V)
\end{eqnarray*}
above equation shows that (a)$\Leftrightarrow$(b).
Further, in view of (3.4) we have
\begin{eqnarray*}
g(\mathcal{T}_V BX,JW)&+&g(\mathcal{A}_{CX} V,JW)\\
&=&g(\mathcal{H}(\nabla_V BX),JW)+g(\mathcal{H}(\nabla_V CX),JW)\\
&=&g(\nabla_V BX,JW)+g(\nabla_V CX,JW)\\
&=&g_B(\pi_*\nabla_V BX,\pi_*JW)+g_B(\pi_*\nabla_V CX,\pi_*JW)\\
&=&-g_B((\nabla\pi_*)(V,BX),\pi_*JW)-g_B((\nabla\pi_*)(V,CX),\pi_*JW)\\
       &+&g_B(\nabla^\pi_V \pi_* CX,\pi_*JW)\\
\end{eqnarray*}
Using (4.1) above equation reduces to
\begin{eqnarray*}
g(\mathcal{T}_V BX,JW)&+&g(\mathcal{A}_{CX} V,JW)\\
&=&-g_B((\nabla\pi_*)(V,BX),\pi_*JW)-g_B((\nabla\pi_*)(V,CX),\pi_*JW)\\
&=&-g_B((\nabla\pi_*)(V,JX),\pi_*JW)
\end{eqnarray*}
and hence we conclude that (b)$\Leftrightarrow$(c).
\hfill$\Box$

Now, if $\pi$ is a Lagrangian Riemannian submersion, then (3.3) implies that $TB=\pi_*(J(ker\pi_*))$. Hence, we have the following

\begin{Corollary}
Let $(M,g,J)$ be a locally conformal Kaehler manifold and $\pi$ be a Lagrangian Riemannian submersion from $M$ onto a Riemannian manifold $(B,g_B)$. Then the following assertions are equivalent to each other:
\begin{description}
  \item[(a)] $(ker\pi_*)$ defines a totally geodesic foliation on $M$
  \item[(b)] $\mathcal{T}_V JW=0$
  \item[(c)] $(\nabla\pi_*)(V,JX)=0$
\end{description}
for $X \in \Gamma((ker\pi_*)^\bot)$ and $V,W \in \Gamma(ker\pi_*).$
\end{Corollary}

{\bf Proof.} $(a)\Leftrightarrow(b)$ is obvius from Theorem 4.3. We prove $(b)\Leftrightarrow(c)$. In view of equation (3.4) and Lemma 2.4 and taking into account the orthogonality of $(ker\pi_*)$ and $(ker\pi_*)^\bot$, we have
\begin{eqnarray*}g(\nabla_V JW,JX)&=&-g(JW,\nabla_V JX)\\
&=&-g_B(\pi_* JW,\pi_* \nabla_V JX)\\
&=&g_B(\pi_* JW,(\nabla\pi_*)(V,JX))\\
g(\mathcal{T}_V JW,JX)&=&g_B(\pi_* JW,(\nabla\pi_*)(V,JX))
\end{eqnarray*}
since $\mathcal{T}_V JW \in \Gamma(ker\pi_*)$, we get $(b)\Leftrightarrow(c)$.
\hfill$\Box$

\begin{Definition} \cite{3} A differential map $\pi$ between two Riemannian manifolds is called totally geodesic if $\nabla\pi_*=0$. \end{Definition}For a Lagrangian Riemannian submersion we have the following characterization result.

\begin{Theorem}
Let $(M,g,J)$ be a locally conformal Kaehler manifold and $\pi$ be a Lagrangian Riemannian submersion from $M$ onto a Riemannian manifold $(B,g_B)$. Then $\pi$ is a totally geodesic map if and only if
\begin{eqnarray*}\mathcal{T}_V JW+\frac{1}{2}\omega(W)JV+\frac{1}{2}g(V,W)A=0\end{eqnarray*}
and
\begin{eqnarray*}\mathcal{A}_X JW+\frac{1}{2}\omega(W)JX+\frac{1}{2}\Omega(X,W)B_1=0\end{eqnarray*}
for $X,Y \in \Gamma((ker\pi_*)^\bot)$ and $V,W \in \Gamma(ker\pi_*).$
\end{Theorem}

{\bf Proof.} We know that second fundamental form of a Riemannian submersion satisfies
\begin{eqnarray}
(\nabla\pi_*)(X,Y)=0   \hspace {1 cm}   \forall  X,Y \in \Gamma((ker\pi_*)^\bot) \end{eqnarray}
In the light of equations (2.1), (3.4) and (4.1), for any $V,W \in (ker\pi_*)$, we have
\begin{eqnarray}
(\nabla\pi_*)(V,W)&=& \nabla^\pi_V \pi_*(W)-\pi_*(\nabla_V W) \nonumber \\
&=&-\pi_*(\nabla_V W) \nonumber \\
&=&\pi_*(J(J\nabla_V W)) \nonumber \\
&=&\pi_*(J(\nabla_V JW+\frac{1}{2}\omega(W)JV+\frac{1}{2}g(V,W)A)) \nonumber \\
&=&\pi_*(J(\mathcal{T}_V JW+\frac{1}{2}\omega(W)JV+\frac{1}{2}g(V,W)A))
\end{eqnarray}
On the other hand, for any $X \in \Gamma((ker\pi_*)^\bot)$ and $W \in (ker\pi_*)$, in view of (3.4), we have
\begin{eqnarray}
(\nabla\pi_*)(X,W)&=& \nabla^\pi_X \pi_*(W)-\pi_*(\nabla_X W) \nonumber \\
&=&-\pi_*(\nabla_X W) \nonumber \\
&=&\pi_*(J(J\nabla_X W)) \nonumber \\
&=&\pi_*(J(\nabla_X JW+\frac{1}{2}\omega(W)JX+\frac{1}{2}\Omega(X,W)B_1)) \nonumber \\
&=&\pi_*(J(\mathcal{A}_X JW+\frac{1}{2}\omega(W)JX+\frac{1}{2}\Omega(X,W)B_1))
\end{eqnarray}
and hence as $J$ is non-singular, result follows from (4.2),(4.3) and (4.4).
\hfill$\Box$

\section{Decomposition Theorems}
\setcounter{equation}{0}
\renewcommand{\theequation}{5.\arabic{equation}}

In this section we obtain decomposition theorems by using the existence of anti-invariant Riemannian submersions. First, we recall the following result from \cite {14}.

\begin{Definition}
\cite {14} Let $g$ be a Riemannian metric tensor on the manifold $N=M\times B$ and assume that the canonical foliations $\mathcal{D}_M$ and $\mathcal{D}_B$ intersect perpendicularly everywhere. Then $g$ is the metric tensor of
\begin{description}
  \item[(i)] a twisted product $M\times_f B$ if and only if $\mathcal{D}_M$ is a totally geodesic foliation and $\mathcal{D}_B$ is a totally umbilical foliation
  \item[(ii)] a warped product $M\times_f B$ if and only if $\mathcal{D}_M$ is a totally geodesic foliation and $\mathcal{D}_B$ is a spherical foliation, that is, it is umbilical and its mean curvature vector field is parallel
  \item[(iii)] a usual product of Riemannian manifolds if and only if $\mathcal{D}_M$ and $\mathcal{D}_B$ are totally geodesic foliations.
\end{description}
\end{Definition}

We have following decomposition theorem for anti-invariant Riemannian submersions which follows from Theorem 4.1 and Theorem 4.3 in terms of second fundamental form of such submersions.

\begin{Theorem}
Let $(M,g,J)$ be a locally conformal Kaehler manifold and $\pi$ be an anti-invariant Riemannian submersion from $M$ onto a Riemannian manifold $(B,g_B)$. Then $M$ is locally product manifold if and only if
\begin{eqnarray*}g_B((\nabla\pi_*)(X,JY),\pi_* JV)&=&-g(CY,J\mathcal{A}_X V)+\frac{1}{2}\omega(V)g(CY,CX)\\
&&-\frac{1}{2}g(BY,B_1)g(X,JV)-\frac{1}{2}g(X,Y)g(B_1,V)
\end{eqnarray*}
and
\begin{eqnarray*}g_B((\nabla\pi_*)(V,JX),\pi_* JW)=0\end{eqnarray*}
for $X,Y \in \Gamma((ker\pi_*)^\bot)$ and $V,W \in \Gamma(ker\pi_*).$
\end{Theorem}

From Corollary 4.2 and Corollary 4.4, we have the following

\begin{Theorem}
Let $(M,g,J)$ be a locally conformal Kaehler manifold and $\pi$ be a Lagrangian Riemannian submersion from $M$ onto a Riemannian manifold $(B,g_B)$. Then $M$ is locally product manifold if and only if
\begin{eqnarray*}g_B((\nabla\pi_*)(X,JY),\pi_* JV)=-\frac{1}{2}g(JY,B_1)g(X,JV)-\frac{1}{2}g(X,Y)g(B_1,V)
\end{eqnarray*}
and
\begin{eqnarray*}\mathcal{T}_V JW=0\end{eqnarray*}
for $X,Y \in \Gamma((ker\pi_*)^\bot)$ and $V,W \in \Gamma(ker\pi_*).$
\end{Theorem}

Now we obtain a decomposition theorem which is related to the notion of twisted product manifold.

\begin{Theorem}
Let $\pi$ be a Lagrangian Riemannian submersion from a locally conformal Kaehler manifold $(M,g,J)$ to a Riemannian manifold $(B,g_B)$. Then $M$ is a locally twisted product manifold of the form $M_{(ker\pi_*)^\bot} \times_f M_{(ker\pi_*)}$ if and only if
\begin{eqnarray*}\mathcal{T}_V JX=-g(X,\mathcal{T}_V V)\|V\|^{-2} JV \end{eqnarray*}
and
\begin{eqnarray*}g_B(\mathcal{A}_X JY,JV)=\frac{1}{2}g(JY,B_1)g(X,JV)+\frac{1}{2}g(X,Y)g(B_1,V)\end{eqnarray*}
for $X,Y \in \Gamma((ker\pi_*)^\bot)$ and $V,W \in \Gamma(ker\pi_*)$, where $M_{(ker\pi_*)^\bot} \times_f M_{(ker\pi_*)}$ are integral manifold of the distributions $(ker\pi_*)^\bot$ and $(ker\pi_*)$.
\end{Theorem}

{\bf Proof.}
Taking into account the orthogonality of $(ker\pi_*)^\bot$ and $(ker\pi_*)$ and using equation (2.1) and Lemma 2.4, , for any $X \in \Gamma((ker\pi_*)^\bot)$ and $V,W \in \Gamma(ker\pi_*)$, we get
\begin{eqnarray*}
g(\nabla_V W,X)&=&-g(\nabla_V X,W) \\
&=&-g(J\nabla_V X,JW) \\
&=&-g(\nabla_V JX,JW) \\
&=&-g(\mathcal{T}_V JX,JW).
\end{eqnarray*}
This implies that $(ker\pi_*)$ is totally umbilical if and only if \begin{eqnarray} \mathcal{T}_V JX=-X(\lambda)JV \end{eqnarray}where $\lambda$ is some function on $M$. Then, from above equation (5.1), we arrive at
\begin{eqnarray} g(-X(\lambda)JV,JV)&=&g(\mathcal{T}_V JX,JV) \nonumber \\
-X(\lambda)\|V\|^2&=&g(\mathcal{T}_V JX,JV) \nonumber \\
&=&g(\nabla_V JX,JV) \nonumber \\
&=&g(J\nabla_V X,JV) \nonumber \\
&=&-g(X,\mathcal{T}_V V)  \nonumber \\
X(\lambda)&=&g(X,\mathcal{T}_V V)\|V\|^{-2}
\end{eqnarray}
where we have used (2.1).
Thus, from (5.1) and (5.2), we conclude that
\begin{eqnarray*}
\mathcal{T}_V JX=-g(X,\mathcal{T}_V V)\|V\|^{-2}JV
\end{eqnarray*}
then result follows from Corollary 4.2.
\hfill$\Box$

We now prove the non existence of a twisted product manifold of the form $M_{(ker\pi_*)^\bot} \times_f M_{(ker\pi_*)}$ for Lagrangian Riemannian submersions.

\begin{Theorem}
Let $(M,g,J)$ be a locally conformal Kaehler manifold  and $(B,g_B)$ be a Riemannian manifold. Then there do not exist Lagrangian Riemannian submersion from $M$ onto $B$ such that $M$ is a locally proper twisted product manifold of the form $M_{(ker\pi_*)^\bot} \times_f M_{(ker\pi_*)}$.
\end{Theorem}

{\bf Proof.}
               Suppose that $\pi$ be a Lagrangian Riemannian submersion from a locally conformal Kaehler manifold $(M,g,J)$ onto a Riemannian manifold $(B,g_B)$ and $M$ is a locally twisted product of the form $M_{(ker\pi_*)^\bot} \times_f M_{(ker\pi_*)}$. Then by definition 5.1 $M_{(ker\pi_*)}$ is a totally geodesic foliation and $M_{(ker\pi_*)^\bot}$ is a totally umbilical foliation. Let us denote the second fundamental form of $M_{(ker\pi_*)^\bot}$ by $h$, then we get $g(\nabla_X Y,V)=g(h(X,Y),V)$ for any $X,Y \in \Gamma((ker\pi_*)^\bot)$ and $V \in \Gamma(ker\pi_*)$. Since $M_{(ker\pi_*)^\bot}$ is totally umbilical foliation then we have
\begin{eqnarray}
g(\nabla_X Y,V)=g(H,V)g(X,Y)
\end{eqnarray}
where $H$ is the mean curvature vector field of $M_{(ker\pi_*)^\bot}$.

On the other hand, taking account of orthogonality of $(ker\pi_*)^\bot$ and $(ker\pi_*)$ and in the light of (2.1) and lemma 2.4, we have
\begin{eqnarray}
g(\nabla_X Y,V)&=&-g(Y,\nabla_X V) \nonumber \\
&=&-g(JY,J\nabla_X V) \nonumber \\
&=&-g(JY,\mathcal{A}_X JV+\frac{1}{2}\omega(V) JX)
\end{eqnarray}
Thus, from (5.3) and (5.4), we get
\begin{eqnarray}
g(H,V)g(X,Y)&=&-g(JY,\mathcal{A}_X JV+\frac{1}{2}\omega(V) JX) \nonumber \\
g(H,V)g(JY,JX)&=&-g(JY,\mathcal{A}_X JV+\frac{1}{2}\omega(V) JX) \nonumber \\
-g(H,V)\|X\|^2&=&g(\mathcal{A}_X JV+\frac{1}{2}\omega(V) JX,JX) \nonumber \\
&=&g(\nabla_X JV+\frac{1}{2}\omega(V) JX,JX) \nonumber \\
&=&g(J\nabla_X V,JX) \nonumber \\
&=&-g(V,\nabla_X X) \nonumber
\end{eqnarray}
So, we have
\begin{eqnarray*}
g(H,V)\|X\|^2=g(V,\mathcal{A}_X X)
\end{eqnarray*}
But, from (2.5) it is known that $\mathcal{A}_X X=0$, which implies $g(H,V)\|X\|^2=0$. Since $g$ is a Riemannian metric and $H \in \Gamma(ker\pi_*)$, we conclude that $H=0$ which proves that $(ker\pi_*)^\bot$ is totally geodesic and hence $M$ is usual product of Riemannian manifolds. This completes the proof.
\hfill $\Box$

\centerline{}
{\bf Acknowledgements.} The author is thankful to Department of Science and Technology, Government of India, for its financial assistance provided through Inspire Fellowship No. DST/INSPIRE Fellowship/2009/[xxv] to carry out this research work.


\begin{thebibliography}{99}

\bibitem{3}{Baird P. and Wood J.C., \em Harmonic morphisms between Riemannian manifolds}, London Mathematical Society Monographs, \textbf{29}, Oxford University Press, The clarendon Press, Oxford, (2003).

\bibitem{6}{Dragomir S. and Ornea, L., \em Locally Conformal Kahler Geometry}, Basel: Birkhauser, (1998).
\bibitem{101}{Gray A., \em Pseudo-Riemannian almost product manifolds and submersions}, J. Math. Mech., {\bf 16} (1967), 715-737.

\bibitem{102}{Ianus S., Mazzocco R. and Vilcu G. E., \em Riemannian submersions from quaternionic manifolds}, Acta Appl. Math., {\bf 104(1)} (2008), 83-89.



\bibitem{4}{Neill B. O'., \em The fundamental equations of a submersion}, Michigan Math. J., {\bf 13}(1996), 459-469.


\bibitem{103}{Park K. S., \em H-slant submersions}, Bull. Korean Math. Soc., {\bf 49(2)} (2012), 329-338.
\bibitem{104}{Park K. S., \em H-semi-invariant submersions}, Taiwanese Journal of Mathematics, {\bf 16(5)} (2012), 1865-1878.


\bibitem{14}{Ponge R. and Reckziegel H., \em Twisted products in pseudo-Riemannian geometry}, Geom. Dedicata, {\bf 48(1)} (1993), 15-25.

\bibitem{11}{Sahin B., \em Anti-invariant Riemannian submersions from almost Hermitian manifolds}, Central European Journal of Mathematics, {\bf 8(3)} (2010), 437-447.



\bibitem{105}{ Sahin B., \em Slant submersions from almost Hermitian manifolds}, Bull.
Math. Soc. Sci. Math. Roumanie Tome, {\bf 54(102) No. 1}, (2011), 93-105.
\bibitem{106}{Sahin B., \em Semi-invariant Riemannian submersions from almost Hermitian manifolds}, Taiwanese math. J., {\bf 17(2)} April 2013, 629-659.



\bibitem{13}{Vaisman I., \em On Locally Conformal Almost Kaehler Manifolds}, Israel J. Math., {\bf 24} (1976), 338-351.





\end{thebibliography}
\end{document}